\documentclass[12pt]{amsart}

\usepackage{amsmath,amsthm,amssymb}
\font\teneufm=eufm10
\font\seveneufm=eufm7
\font\fiveeufm=eufm5
\newfam\frakturfam

\textfont\frakturfam=\teneufm
\scriptfont\frakturfam=\seveneufm
\scriptscriptfont\frakturfam=\fiveeufm

\newtheorem{pr}{Proposition}

\newtheorem{lm}{Lemma}
\newtheorem{theor}{Theorem}
\newtheorem{co}{Corollary}

\def\bee{\begin{eqnarray}}
\def\bes{\begin{eqnarray*}}
\def\eee{\end{eqnarray}}
\def\ees{\end{eqnarray*}}

\def\a{\alpha}

\def\Proof{{\sl Proof.}\ }

\title{Free Poisson fields and their automorphisms}

\begin{document}
\date{}
\maketitle

\begin{center}

{\bf Leonid Makar-Limanov}\footnote{Supported
by an NSA grant H98230-09-1-0008, by an NSF grant DMS-0904713, and a Fulbright fellowship awarded by the United States--Israel Educational Foundation; The Weizmann Institute of Science, Rehovot, Israel, University of Michigan, Ann Arbor, and
 Wayne State University, Detroit, MI 48202, USA,
e-mail: {\em lml@math.wayne.edu}}
and
{\bf Ualbai Umirbaev}\footnote{Supported by an NSF grant DMS-0904713 and by a grant of Kazakhstan; Eurasian National University,
 Astana, Kazakhstan and
 Wayne State University,
Detroit, MI 48202, USA,
e-mail: {\em umirbaev@math.wayne.edu}}

\end{center}

\begin{abstract}
Let $k$ be an arbitrary field of characteristic $0$. We prove that the group of automorphisms of a free Poisson field $P(x,y)$ in two variables $x,y$ over $k$ is isomorphic to the Cremona group $\mathrm{Cr}_2(k)$. We also prove that the universal enveloping algebra $P(x_1,\ldots,x_n)^e$
of a free Poisson field $P(x_1,\ldots,x_n)$ is a free ideal ring and give a characterization of the Poisson dependence of two elements of $P(x_1,\ldots,x_n)$ via universal derivatives.
\end{abstract}


\noindent


\section{Introduction}

\hspace*{\parindent}

It is well known \cite{Czer,Jung,Kulk,Makar} that the automorphisms
of polynomial algebras and free associative algebras in two
variables are tame. It is also known \cite{Umi25,Um19} that
polynomial algebras and free associative algebras in three
variables in the case of characteristic zero have wild
automorphisms. It was recently proved \cite{MLTU} that the
automorphisms of free Poisson algebras in two variables  are tame in characteristic zero. Note that the Nagata
automorphism \cite{Nagata,Umi25} gives an example of a wild
automorphism of a free Poisson algebra in three variables.

The famous Max Noether Theorem on the Cremona groups (see \cite{Nagata3,AS}) says that the Cremona group $Cr_2(k)$ of birational automorphisms of the projective plane $\mathbb{P}_k^2$ over an algebraically closed field $k$ is generated by the standard quadratic transformation and the projective transformations. The Cremona group $Cr_2(k)$ is isomorphic to the group of automorphisms of the field of rational functions $k(x,y)$ over $k$.
There are some generalizations of Noether's theorem to the case in which the ground field $k$ is not algebraically closed (see, for example \cite{Manin}). Defining relations of the Cremona group $Cr_2(k)$ are described in \cite{Gizat} (see also \cite{Iskovskikh}).

Let $P\langle x_1,\ldots,x_n\rangle$ be a free Poisson algebra over a field $k$ freely generated by $x_1,\ldots,x_n$  \cite{Shest}. Denote by $P(x_1,\ldots,x_n)$ the field of fractions of the commutative and associative algebra  $P\langle x_1,\ldots,x_n\rangle$. The Poisson bracket $\{\cdot,\cdot\}$ on $P\langle x_1,\ldots,x_n\rangle$ can be uniquely extended to a Poisson bracket on $P(x_1,\ldots,x_n)$. The field $P(x_1,\ldots,x_n)$ with this Poisson bracket is called the free Poisson field in the variables $x_1,\ldots,x_n$ over $k$ \cite{MLSh}.

The main goal of this paper is to prove that the group of automorphisms of the free Poisson field $P(x,y)$ over a filed $k$ of characteristic $0$ is isomorphic to the Cremona group $Cr_2(k)$. Notice that there are several different proofs \cite{MLTU,MLU,MLSh,UU10} of the tameness of automorphisms of the free Poisson algebra $P\langle x,y\rangle$ in characteristic zero. In \cite{UU10}  the technic of universal derivations or Fox derivatives is used. We adopt this method to the case of free fields.

For any Poisson algebra $P$ denote by $P^e$ its universal enveloping algebra.
The main property of the universal   enveloping algebra $P^e$ is that the notion of a Poisson $P$-module (see, for example \cite{Oh}) is equivalent to the notion of a left module over the associative algebra $P^e$ (see, for example \cite{Jacobson}). The universal enveloping algebras of Poisson algebras were first studied in \cite{Oh}. A linear basis of the universal enveloping algebra is constructed in \cite{OPC}  for Poisson polynomial algebras. The structure of the universal enveloping algebras of free Poisson algebras were studied in \cite{UU10}.

It is well known that the universal enveloping algebras of free Lie algebras are free associative algebras. P.\,Cohn  \cite{Cohn} proved that every left ideal of a free associative algebra is a free left module, i.e., a free associative algebra is a free ideal ring (also fir or $\mathrm{FI}$-ring \cite{Cohn}). In fact, it follows from this theorem \cite{Um} that subalgebras of free Lie algebras are free \cite{Shir1,Witt} and automorphisms of finitely generated free Lie algebras are tame \cite{Cohn2}. Unfortunately the universal enveloping algebras of free Poisson algebras are not free ideal rings \cite{UU10}. In this paper we prove that the universal enveloping algebras of free Poisson fields  are free ideal rings.

Recall that any two elements of a free Poisson field are Poisson dependent if and only if they are algebraically dependent \cite{MLSh}. In \cite{UU10} a characterization of Poisson dependence of two elements of free Poisson algebras via universal derivations is given. We extend this result to free Poisson fields.

This paper is organized as follows. In Section 2 we recall the definition of the universal enveloping algebras of Poisson algebras and construction of free Poisson algebras and free Poisson fields. In Section 3 we describe the structure of the universal enveloping algebras of free Poisson fields and prove that they are free ideal rings. Section 4 is devoted to the study of Poisson dependence of two elements. Developed technics are used to describe automorphisms of free Poisson fields in two variables.

\section{Definitions and notations}

\hspace*{\parindent}

Recall that a vector space $P$ over $k$ endowed with two bilinear
operations $x\cdot y$ (a multiplication) and $\{x,y\}$ (a Poisson
bracket) is called {\em a Poisson algebra} if $P$ is a commutative
associative algebra under $x\cdot y$, $P$ is a Lie algebra under
$\{x,y\}$, and $P$ satisfies the following identity (the Leibniz
identity): \bes \{x, y\cdot z\}=\{x,y\}\cdot z + y\cdot \{x,z\}.
\ees

Let $P$ be a Poisson algebra over $k$.
A vector space $V$ over $k$ is called a {\em Poisson module} over $P$
(or {\em Poisson $P$-module}) if there are two bilinear maps
\bes
P\,\times \,V\longrightarrow V, ((x,v)\mapsto x\cdot v), \ \ \ \ P\,\times \,V\longrightarrow V, ((x,v)\mapsto \{x,v\}),
\ees
such that the relations
\bes
(x\cdot y)\cdot v=x\cdot (y\cdot v),\\
\{\{x,y\},v\}=\{x,\{y,v\}\}-\{y,\{x,v\}\},\\
\{x\cdot y,v\}=y\cdot \{x,v\}+x\cdot \{y,v\},\\
\{x,y\}\cdot v=\{x,y\cdot v\}-y\cdot \{x,v\}
\ees
hold for all $x,y\in P$ and $v\in V$. A Poisson module $V$ is called unitary if $1\cdot v=v$ for all $v\in V$.

The universal enveloping algebra $P^e$ of $P$ is an associative algebra for which
the category of Poisson modules over $P$ and the category of (left) modules over $P^e$ are equivalent.
It can be
given by generators and defining relations.  Let $m_P=\{m_a | a\in P\}$ and $h_P=\{h_a | a\in P\}$ be two copies of the vector space $P$ endowed with two linear isomorphisms $m : P \longrightarrow m_P \ (a\mapsto m_a)$ and $h : P \longrightarrow h_P \ (a\mapsto h_a)$. Then $P^e$ is an associative algebra over $k$, with an identity $1$, generated by two linear spaces $m_P$ and $h_P$ and defined by the relations
\bes
m_{x y}=m_x m_y, \\
h_{\{x,y\}}=h_x h_y - h_y h_x,\\
h_{xy}=m_y h_x + m_x h_y,\\
 m_{\{x,y\}}=h_x m_y-m_y h_x
\ees
for all $x,y\in P$ \cite{UU10}. Additionally we put $m_1=1$ since we consider only unitary Poisson modules.

If $V$ is an arbitrary Poisson $P$-module, then $V$ becomes a left $P^e$-module under the actions
\bes
m_xv=x\cdot v, \  h_x v= \{x,v\},
\ees
for all $x\in P$ and $v\in V$.
Conversely, if $V$ is a left $P^e$-module then the same formulae turn $V$ to a Poisson $P$-module.

The first example of a Poisson $P$-module is $V=P$ under the actions $x\cdot v$ and $\{x,v\}$ which provides us with a representation $\rho$ of $P^e$ by homomorphisms of $V$.
It is clear that $\rho(m_x)(v)$ is $xv$ and hence we have an imbedding
\bes
m : P\longrightarrow P^e \ \ \ \ \ (x\mapsto m_x).
\ees
Therefore we can identify $m_x$ with $x$. After this identification the defining relations of $P_e$ are
\bee\label{f1}
h_{\{x,y\}}=h_x h_y - h_y h_x,
\eee
\bee\label{f2}
h_{xy}=y h_x + x h_y,
\eee
\bee\label{f3}
\{x,y\}=h_x y-y h_x,
\eee
for all $x,y\in P$. It follows from (\ref{f3}) that
\bee\label{f4}
h_x y-y h_x = x h_y- h_y x.
\eee
for all $x,y\in P$. From (\ref{f2}) we obtain that $h_1=0$ and if $x\in P$ is invertible then
\bee\label{f5}
h_{x^{-1}}=-x^{-2}h_x.
\eee

A word of caution, $\rho$ is not necessarily an exact representation. A question for which $P \ \ \rho$ is exact seems to be an interesting problem.

Let $g$ be a Lie algebra with a linear
basis $e_1,e_2,\ldots,e_k,\ldots$. The Poisson symmetric algebra $PS(g)$ of
$g$ is the usual polynomial algebra $k[e_1,
e_2,\ldots,e_k,\ldots]$ endowed with the Poisson bracket defined
by  \bes \{e_i,e_j\}=[e_i,e_j] \ees for all $i,j$, where $[x,y]$
is the multiplication in the Lie algebra $g$.

Denote by $P\langle x_1,x_2,\ldots,x_n\rangle$ the free Poisson algebra over $k$ in the variables $x_1,x_2,\ldots,x_n$.
From now on let $g=\mathrm{Lie}\langle x_1,x_2,\ldots,x_n\rangle$ be the free Lie algebra with free (Lie)
generators $x_1,x_2,\ldots,x_n$. It is well-known (see, for
example \cite{Shest}) that the Poisson symmetric algebra $PS(g)$ is the free Poisson algebra $P\langle x_1,x_2,\ldots,x_n\rangle$ in the variables $x_1,x_2,\ldots,x_n$.

Let us choose a multihomogeneous linear basis \bes
x_1,x_2,\ldots,x_n,\,[x_1,x_2],\ldots,[x_1,x_n],\ldots,[x_{n-1},x_n],\,[[x_1,x_2],x_3],\ldots
\ees of a free Lie algebra $g$ and denote the elements of this
basis by
\bee
\label{f6} e_1, e_2, \ldots, e_m, \ldots.
\eee

The algebra $P\langle x_1,x_2,\ldots,x_n\rangle$ coincides with
the polynomial algebra on the elements (\ref{f6}).
Consequently, the set of all words of the form
\bee\label{f7}
u=e^{\a}=e_1^{i_1}e_2^{i_2}\ldots e_m^{i_m},
\eee
where $0\leq i_k$, $1\leq k\leq m$, and $m\geq 0$, forms a linear basis of $P$.

By $\deg$ we denote the
Poisson degree function on $P\langle x_1,x_2,\ldots,x_n\rangle$ as the polynomial algebra on the elements (\ref{f6}), i.e.
$\deg(e_i)$ is the Lie degree of this element. Note that
\begin{eqnarray*}
\deg\,\{f,g\}= \deg\,f+\deg\,g
\end{eqnarray*}
if $f$ and $g$ are homogeneous and $\{f,g\}\neq 0$.
By
$\deg_{x_i}$ we denote the degree function on $P$ with respect to
$x_i$. i.e.
$\deg_{x_i}(e_j)$ is the number of appearances of $x_i$ in the Lie monomial $e_j$.
If $f$ is homogeneous with respect to each $\deg_{x_i}$, where
$1\leq i \leq n$, then $f$ is called multihomogeneous.
The basis (\ref{f7}) is multihomogeneous since so is (\ref{f6}).

Denote by $P(x_1,x_2,\ldots,x_n)$ the field of fractions of the polynomial algebra
$k[e_1,\ldots,e_n,\ldots]$ in the variables (\ref{f6}). The Poisson bracket $\{\cdot,\cdot\}$ on $k[e_1,\ldots,e_n,\ldots]=P\langle x_1,x_2,\ldots,x_n\rangle$ can be uniquely extended to a Poisson bracket on
$P(x_1,x_2,\ldots,x_n)$  and
\bes
\{\frac{a}{b},\frac{c}{d}\}=\frac{\{a,c\}bd-\{a,d\}bc-\{b,c\}ad+\{b,d\}ac}{b^2d^2}
\ees
for all $a,b,c,d\in P\langle x_1,x_2,\ldots,x_n\rangle$ with $bd\neq 0$.

The field $P(x_1,x_2,\ldots,x_n)$
with this Poisson bracket is called the {\em free Poisson field} over $k$ in the variables $x_1,x_2,\ldots,x_n$ \cite{MLSh}.

\section{Enveloping algebras of free Poisson fields}

\hspace*{\parindent}

\begin{pr}\label{p1}
Let $Q=P(x_1,x_2,\ldots,x_n)$ be the free Poisson field over $k$ in the variables $x_1,x_2,\ldots,x_n$ and let $Q^e$ be its universal enveloping algebra. Then the following statements are true:

(i) The subalgebra $A$ (with identity) of $Q^e$ generated by $h_{x_1},h_{x_2},\ldots, h_{x_n}$ is a free associative algebra freely generated by $h_{x_1},h_{x_2},\ldots, h_{x_n}$;

(ii) The left vector space $Q^e$ over $Q$ is isomorphic to the left vector space $Q\otimes_k A$ over $Q$.
\end{pr}
\Proof
Let $M = P(x_1,x_2,\ldots,x_n,y)$ be the free Poisson field in the variables $x_1,x_2,\ldots,x_n,y$.
Then $M$ is a Poisson $Q$-module under the actions $x\cdot v$ and $\{x,v\}$ and correspondingly a $Q^e$ module.
Denote by $N$ the submodule of $M$ generated by $y$.

Consider the elements
\bee\label{f7'}
h_{x_{i_1}}h_{x_{i_2}}\ldots h_{x_{i_k}}(y) = \{x_{i_1},\{x_{i_2},\ldots,\{x_{i_k},y\}\ldots\}\}
\eee
of $N$.
They are a subset of the multihomogeneous liner basis of $Lie\langle x_1,x_2,\ldots,x_n,y\rangle$ and are linearly independent.
Consequently, the elements of the form
\bee\label{f8}
h_{x_{i_1}}h_{x_{i_2}}\ldots h_{x_{i_k}}
\eee
are linearly independent in $Q^e$ and $A$ is a free associative algebra freely generated by $h_{x_1},h_{x_2},\ldots, h_{x_n}$.

The relations (\ref{f3}) allow to express every element of $Q^e$ as a linear combination
$\sum_i q_i w_i$, where $q_i\in Q$ and $w_i$ are different elements of the form (\ref{f8}). Then $\sum_i q_i w_i(y) \neq 0$. Indeed, $M$ is a polynomial ring and $q_i$ are polynomials in variables which are algebraically independent with the variables (\ref{f7'}) since all these variables are different elements of the liner basis of $Lie\langle x_1,x_2,\ldots,x_n,y\rangle$.
This means that (ii) is true.  $\Box$

\begin{co}\label{c1}
Every nonzero element $u$ of the universal enveloping algebra $Q^e$ can be uniquely written in the form
\bee\label{f9}
u=\sum_{i=1}^k q_i w_i,
\eee
where $0\neq q_i\in Q$ for all $i$ and $w_1,w_2,\ldots,w_k$ are different elements of the form (\ref{f8}).
\end{co}

Let $u$ be an element of $Q^e$ written in the form (\ref{f9}). Put $\mathrm{hdeg}\,u= max_{i=1}^k \mathrm{hdeg}\,w_i$ where $\mathrm{hdeg}$ is the homogeneous degree on the free algebra $A$ and $\mathrm{hdeg}\,0=-\infty$.  We say that $u$ is homogeneous with respect to $\mathrm{hdeg}$ if $\mathrm{hdeg}\,w_1=\mathrm{hdeg}\,w_2=\ldots =\mathrm{hdeg}\,w_k$.

Put $h_{x_i}<h_{x_j}$ if $i<j$. Let $u,v$ be two elements of the form (\ref{f8}). Then put $u<v$ if $\mathrm{hdeg}\,u<\mathrm{hdeg}\,v$ or $\mathrm{hdeg}\,u=\mathrm{hdeg}\,v$ and $u$ precedes $v$ in the lexicographical order.

Let $u$ be an element of the form (\ref{f9}). We may assume that $w_1<w_2<\ldots <w_k$. Then $w_k$ is called the {\em leading monomial} of $u$ and $q_k$ is called the {\em leading coefficient} of $u$. We write $w_k=\mathrm{ldm}(u)$ and $q_k=\mathrm{ldc}(u)$. The {\em leading term} of $u$ is defined to be $\mathrm{ldt}(u)=\mathrm{ldc}(u) \mathrm{ldm}(u)$.

\begin{lm}\label{l1}
If $u$ and $v$ are arbitrary nonzero elements of $Q^e$ then
\bes
\mathrm{ldc}(uv)=\mathrm{ldc}(u)\mathrm{ldc}(v) \ \ \ {\it and} \ \ \  \mathrm{ldm}(uv)=\mathrm{ldm}(u)\mathrm{ldm}(v).
\ees
\end{lm}
\Proof
Note that if $u$ and $v$ are two elements of the form (\ref{f9}) then to put the product $uv$ into the form (\ref{f9}) again we need to use only the relations of type (\ref{f3}) which imply that $h_{x_i}$ and $y\in P$ commute modulo terms of smaller degrees in the variables $h_{x_1},h_{x_2},\ldots,h_{x_n}$. Consequently, we can put $uv$ into the form (\ref{f9}) with the leading monomial $\mathrm{ldm}(u)\mathrm{ldm}(v)$ and the leading coefficient $\mathrm{ldc}(u)\mathrm{ldc}(v)$.
 $\Box$

It follows directly from Lemma \ref{l1} that
\bes
\mathrm{hdeg}\,uv=\mathrm{hdeg}\,u+\mathrm{hdeg}\,v
\ees
for every $u$ and $v$ from $Q^e$, i.e., $\mathrm{hdeg}$ is a degree function on $Q^e$.

Denote by $U_i$ the subset of all elements $u$ of $Q^e$ with $\mathrm{hdeg}\,u\leq i$. Then
\bes
Q=U_0\subset U_1\subset U_2\subset\ldots \subset U_k\subset \ldots ,
\ees
is a filtration of $Q^e$, i.e., $U_i U_j\subseteq U_{i+j}$ for all $i,j\geq 0$ and $\cup_{i\geq 0} U_i= Q^e$. Put
\bes
\mathrm{gr}\,Q^e= \mathrm{gr}\,U_0\oplus \mathrm{gr}\,U_1\oplus \mathrm{gr}\,U_2\oplus\ldots \oplus \mathrm{gr}\,U_k\oplus\ldots,
\ees
where $\mathrm{gr}\,U_0=Q$ and $\mathrm{gr}\,U_i=U_i/U_{i-1}$ for all $i\geq 1$.
Denote by $\varphi_i : U_i\rightarrow \mathrm{gr}\,U_i$ the natural projection for every $i\geq 1$ and define
\bee\label{f10}
\varphi=\{\varphi_i\}_{i\geq 0} : Q^e\rightarrow \mathrm{gr}\,Q^e
\eee
by $\varphi(u)=\varphi_i(u)$ if $u\in U_i\setminus U_{i-1}$ for every $i\geq 1$ and $\varphi(u)=u$ if $u\in Q$.

The multiplication on $Q^e$ induces a multiplication on $\mathrm{gr}\,Q^e$ and the graded vector space $\mathrm{gr}\,Q^e$ becomes an algebra.

Recall that $A$ is the free associative subalgebra of $Q^e$ generated by $h_{x_1},h_{x_2},\ldots,h_{x_n}$ (Proposition \ref{p1}). Let $B=Q\otimes_k A$ be the tensor product of associative algebras over $k$. Then $B$ is a free associative algebra over $Q$ freely generated by $h_{x_1},h_{x_2},\ldots,h_{x_n}$.

\begin{pr}\label{p2}
The graded algebra $\mathrm{gr}\,Q^e$ is isomorphic to $B=Q\otimes_k A$.
\end{pr}
\Proof By (\ref{f3}), $Q$ is in the center of the algebra $\mathrm{gr}\,Q^e$ and $\mathrm{gr}\,Q^e$ is generated by $\varphi(h_{x_1}),\varphi(h_{x_2}),\ldots,\varphi(h_{x_n})$ as an algebra over $Q$. Note that $B=Q\otimes_k A$ is a free associative algebra over $Q$. Hence there is a $Q$-algebra homomorphism $\psi : B \rightarrow \mathrm{gr}\,Q^e$ such that $\psi(h_{x_i})=\varphi(h_{x_i})$ for all $i$.

Let $T_s$ be the space of $\mathrm{hdeg}$ homogeneous elements of $Q^e$ of degree $s\geq 0$ and $B_s$ be the space of homogeneous elements of degree $s$ of $B$. There is an obvious isomorphism between the $Q$-linear spaces $T_s$ and $B_s$ established by Corollary \ref{c1}. Note that $U_s=U_{s-1}+T_s$, $U_s/U_{s-1}\simeq T_s\simeq B_s$, and $\psi_{|B_s} : B_s\rightarrow \mathrm{gr}\,U_s$ is an isomorphism of $Q$-modules.  Consequently, $Ker\,(\psi)=0$ and $\psi$ is an isomorphism of algebras. $\Box$

Therefore $\mathrm{gr}\,Q^e$ is a free associative algebra over $Q$ and $\mathrm{gr}\,Q^e$ is a free ideal ring \cite{Cohn}, i.e., every left (right) ideal of $\mathrm{gr}\,Q^e$ is a free left (right) $\mathrm{gr}\,Q^e$-module of unique rank. From now on we will identify $\mathrm{gr}\,Q^e$ with $B=Q\otimes_k A$ by means of the isomorphism from Proposition \ref{p2}. Obviously, in this identification the elements $\varphi(h_{x_1}),\ldots,\varphi(h_{x_n})\in \mathrm{gr}\,Q^e$ correspond to $h_{x_1},\ldots,h_{x_n}\in B$. Denote by $\mathrm{Deg}$ a $Q$-degree function on $B$ such that
 $\mathrm{Deg}(h_{x_i})=1$ for all $i$. Notice that for any $u\in Q^e$ we have $\mathrm{hdeg}(u)=\mathrm{Deg}(\varphi(u))$.

\begin{theor}\label{t1} Let $Q=P(x_1,x_2,\ldots,x_n)$ be the free Poisson field over $k$ in the variables
$x_1,x_2,\ldots,x_n$  and $Q^e$ be its universal enveloping algebra. Then $Q^e$ satisfies the weak algorithm for $\mathrm{hdeg}$ and is a free ideal ring.
\end{theor}
\Proof We consider only the left dependence in $Q^e$ since the right dependence can be treated similarly.
It is sufficient to prove that $Q^e$ satisfies the weak algorithm for $\mathrm{hdeg}$ \cite{Cohn}, i.e., for any finite set of left $Q^e$-dependent elements $s_1,s_2,\ldots,s_k$ of $Q^e$ with
$\mathrm{hdeg}(s_1)\leq \mathrm{hdeg}(s_2)\leq\ldots\leq \mathrm{hdeg}(s_k)$ there exist $i, 1\leq i\leq k$ and $v_1,\ldots,v_{i-1}\in Q^e$ such that $\mathrm{hdeg}(s_i-v_1s_1-\ldots-v_{i-1}s_{i-1})<\mathrm{hdeg}(s_i)$ and
$\mathrm{hdeg}(v_j)+\mathrm{hdeg}(s_j)\leq\mathrm{hdeg}(s_i)$ for all $j, 1\leq j\leq i-1$.

Suppose that
\bes
\sum_{r=1}^k u_r s_r=0.
\ees
Put $m=\max\{\mathrm{hdeg}(u_i)+\mathrm{hdeg}(s_i) | 1\leq i\leq k\}$. Let $j_1,\ldots,j_r$ be the set of indexes $i$ with $\mathrm{hdeg}(u_i)+\mathrm{hdeg}(s_i)=m$. Then the last equality induces
\bes
\varphi(u_{i_1})\varphi(s_{i_1})+\varphi(u_{i_2})\varphi(s_{i_2})+\ldots+\varphi(u_{i_r})\varphi(s_{i_r})=0
\ees
in the free associative algebra $B$ over $Q$. This is a nontrivial left dependence of the homogeneous elements
$\varphi(s_{i_1}),\varphi(s_{i_2}),\ldots,\varphi(s_{i_r})$ in $B$. Since the free associative algebra $B$ over the field $Q$ satisfies the weak algorithm for $\mathrm{Deg}$ \cite{Cohn} it follows that there exist $t$ and $v_{i_1},\ldots,v_{i_{t-1}}\in Q^e$ such that
\bes
\varphi(s_{i_t})= \varphi(v_{i_1})\varphi(s_{i_1})+\ldots+\varphi(v_{i_{t-1}})\varphi(s_{i_{t-1}})
\ees
and $\mathrm{Deg}(\varphi(v_{i_j}))+\mathrm{Deg}(\varphi(s_{i_j}))= \mathrm{Deg}(\varphi(s_{i_t}))$ for all $j, 1\leq j\leq t-1$. Then
\bes
\mathrm{hdeg}(s_{i_t}- v_{i_1}s_{i_1}-\ldots-v_{i_{t-1}}s_{i_{t-1}})<
\mathrm{hdeg}(s_{i_t})
\ees
and $\mathrm{hdeg}(v_{i_j})+\mathrm{hdeg}(s_{i_j})\leq \mathrm{hdeg}(s_{i_t})$ for all $j, 1\leq j\leq t-1$.
$\Box$

The proof of this theorem is constructive. Hence the standard algorithms (see, for example \cite{Um6,Um11,Um}) give the next result.
\begin{co}\label{c3} $(i)$ The left ideal membership problem for $Q^e$ is algorithmically decidable;

$(ii)$ The left dependence of a finite system of elements of $Q^e$ is algorithmically recognizable.
\end{co}

\section{Poisson dependence of two elements}

\hspace*{\parindent}

Let $P$ be an arbitrary Poisson algebra over $k$. Any finite set of elements $p_1,p_2,\ldots,p_k$ of $P$ is called
{\em Poisson dependent} if there exists a nonzero element $f$ of the free Poisson algebra $P\langle x_1,x_2,\ldots,x_k\rangle$ such that $f(p_1,p_2,\ldots,p_k)=0$. Otherwise $p_1,p_2,\ldots,p_k$
are called {\em Poisson free} or {\em Poisson independent}. If $p_1,p_2,\ldots,p_k$ are Poisson free then the Poisson subalgebra of $P$ generated by these elements is a free Poisson algebra in these variables. Similarly we can define Lie dependence and associative dependence \cite{Um6,UU93,Um11,Um}.

Denote by $\Omega_Q$ the left ideal of $Q^e$ generated by $h_{x_1},h_{x_2},\ldots, h_{x_n}$. By Proposition  \ref{p1},
\bes
\Omega_Q = Q^eh_{x_1}\oplus Q^eh_{x_2}\oplus\ldots \oplus Q^eh_{x_n},
\ees
i.e., $\Omega_Q$ is a free left $Q^e$-module. Note that
\bes
Q^e=Q \oplus \Omega_Q.
\ees

Consider
\bes
h : Q \rightarrow \Omega_Q
\ees
such that $h(q)=h_q$ for all $q\in Q$.
It follows from (\ref{f1}) and (\ref{f2}) that $h$ is a derivation of the Poisson algebra $Q$ with coefficients in the Poisson $Q$-module $\Omega_Q$. It is proved in \cite{UU10} that $h$ is the universal derivation of the Poisson algebra $Q$ and $\Omega_Q$ is the universal differential module of $Q$.
It means (see, for example \cite{Um11,Um}) that for any Poisson $Q$-module $V$ and for any derivation $d : Q \rightarrow V$ there exists a unique homomorphism $\tau : \Omega_Q \rightarrow V$ of Poisson $Q$-modules such that $d=\tau h$.

The proof of the next lemma is standard \cite{Um11,Um,UU10}.

\begin{lm}\label{l3} Let $f_1,f_2,\ldots,f_k$ be arbitrary elements of the free Poisson field $Q$.
If $f_1,f_2,\ldots,f_k$ are Poisson dependent then
$h_{f_1},h_{f_2},\ldots,h_{f_k}$ are left dependent over $Q^e$.
\end{lm}
\Proof We repeat the proof given in \cite{UU10} for a free Poisson algebra.
Let $f=f(z_1,z_2,\ldots,z_k)$ be a nonzero element of $T=P\langle z_1,z_2,\ldots,z_k\rangle$ such that
$f(f_1,f_2,\ldots,f_k)=0$ with the minimal degree possible. Then $h_f$ can be written in $\Omega_T$ as
\bes
h_f=u_1 h_{z_1}+u_2 h_{z_2}+\ldots + u_k h_{z_k}.
\ees
and
\bes
0=h_{f(f_1,f_2,\ldots,f_k)}= u_1' h_{f_1}+u_2' h_{f_2}+\ldots + u_k' h_{f_k},
\ees
where $u_i'=u_i(f_1,f_2,\ldots,f_k)$ for all $i$.

Since $f \in T\setminus k$ and hence $f$ is not in the Poisson center of $T$ it is clear that $h_f$ is not identically zero. Therefore we may assume that $u_1=u_1(z_1,z_2,\ldots,z_k)\neq 0$.
If $u_1'\neq 0$ then the last equation gives a nontrivial dependence of $h_{f_1},h_{f_2},\ldots,h_{f_k}$.
Suppose that $u_1'= 0$. Note that $u_1=t+w$, where $t\in T$ and $w\in \Omega_T$, since $U(T)=T\oplus \Omega_T$.
Obviously, $t(f_1,f_2,\ldots,f_k)\in Q$ and it follows from (\ref{f1})--(\ref{f2}) that $w(f_1,f_2,\ldots,f_k)\in \Omega_Q$. Then, $t(f_1,f_2,\ldots,f_k)=0$ and $w(f_1,f_2,\ldots,f_k)=0$ since
$0=u_1'=t(f_1,f_2,\ldots,f_k)+w(f_1,f_2,\ldots,f_k)\in Q\oplus \Omega_Q$.
For any monomial $M \in T$ we have $h_M = \sum_i t_i w_i$ where $t_i \in T, \ w_i$ belong to a free associative algebra generated by $h_{z_1},h_{z_2},\ldots,h_{z_k}$ and $\deg t_i + \mathrm{hdeg}\, t_i = \deg\, M$. Therefore $\deg\, t \leq \deg\, f \ - 1$ and $t = 0$ since the degree of $f$ was minimal possible.
Hence $w\neq 0$ and $w' = w(f_1,f_2,\ldots,f_k)$ provides a "smaller" dependence between $h_{f_1},h_{f_2},\ldots,h_{f_k}$ and after a finite number of steps we get a nontrivial dependence between $h_{f_1},h_{f_2},\ldots,h_{f_k}$ over $Q^e$.
$\Box$

 For every $i\geq 1$ denote by $\partial_{e_i}$ the usual partial derivation of the field of rational functions $Q$ in the variables (\ref{f6}). The next easy lemma is useful.
\begin{lm}\label{l4} $(i)$ For every $q\in Q$ we have
\bes
h_q= \sum_{i\geq 1} \partial_{e_i}(q) h_{e_i};
\ees

$(ii)$ For every $l=l(x_1,x_2,\ldots,x_n)\in g=Lie\langle x_1,x_2,\ldots,x_n\rangle$ we have
\bes
h_l=l(h_{x_1},h_{x_2},\ldots,h_{x_n}) \in L = Lie\langle h_{x_1},h_{x_2},\ldots,h_{x_n}\rangle.
\ees
\end{lm}
\Proof As mentioned above, by (\ref{f1}) and (\ref{f2}), the mapping $h$ is a derivation with respect to both operations on $Q$.
$\Box$
\begin{co}\label{c4} For every $q\in Q$ we have $h_q\in QL$ and $\mathrm{hdeg}(h_q)\leq 1$ if and only if $q\in k(x_1, x_2, \ldots, x_n)$.
\end{co}

We need the next lemma.
\begin{lm}\label{l2}\cite{Um6,Um9,Um11}  Let $H$ be a free Lie algebra over a field $k$ and $U(H)$ be its universal enveloping algebra. Then arbitrary elements $f_1,\ldots,f_k$ of $H$ are Lie dependent if and only if they are left dependent over $U(H)$.
\end{lm}

Recall that $A$ is a free associative algebra freely generated by $h_{x_1},h_{x_2},\ldots,h_{x_n}$.
Hence $L$ is the Lie subalgebra of the Lie algebra $A^{(-)}$. Notice that $B=Q\otimes_k A$ is a free associative algebra over $Q$ and $Q\otimes_k L$ is a free Lie algebra over $Q$ freely generated by $h_{x_1},h_{x_2},\ldots,h_{x_n}$. Furthermore $B$ is the universal enveloping algebra of the Lie algebra $Q\otimes_k L$ over $Q$.

\begin{theor}\label{t2} Let $Q=P(x_1,x_2,\ldots,x_n)$ be the Poisson free field in the variables $x_1,x_2,\ldots,x_n$ over a field $k$ of characteristic zero and $f$ and $g$ be arbitrary elements of $Q$. Then the following conditions are equivalent:

$(i)$ $f$ and $g$ are Poisson dependent;

$(ii)$ $h_f$ and $h_g$ are left dependent over $Q^e$;

$(iii)$ $f$ and $g$ are polynomially dependent, i.e., they are algebraically dependent in the commutative algebra $Q$ over $k$;

$(iv)$ there exists $a\in Q$ such that $f,g\in k(a)$;

$(v)$ $\{f,g\}=0$ in $Q$.
\end{theor}
\Proof The implication $(iii)\rightarrow (iv)$ follows from the generalized L\"uroth Theorem \cite{Sam}. The implications $(iv)\rightarrow (v)$ and $(v)\rightarrow (i)$ are obvious. By Lemma \ref{l3}, (i) implies (ii).
To finish the proof  it is sufficient to show that $(ii)$ implies $(iii)$.

Notice that $h_f=0$ if and only if $f\in k$ by Lemma \ref{l4}. Suppose that $h_f$ and $h_g$ are both nonzero and
\bes
u h_f+ v h_g=0, \ \ (u,v)\neq (0, 0), \  \ u,v\in Q^e.
\ees
It follows that
\bes
\mathrm{hdeg}(u)+\mathrm{hdeg}(h_f)=\mathrm{hdeg}(v)+\mathrm{hdeg}(h_g)
\ees
by Lemma \ref{l1} and the mapping $\varphi$ defined in (\ref{f10}) gives
\bes
\varphi(u) \varphi(h_f)+ \varphi(v) \varphi(h_g)=0, \ \ (\varphi(u),\varphi(v))\neq (0, 0)
\ees
in $B=Q\otimes_k A$. So, $\varphi(h_f)$ and $\varphi(h_g)$ are left dependent over $B$.

By Lemma \ref{l4}, $h_f,h_g\in QL$.
Consequently, $\varphi(h_f), \varphi(h_g)\in Q\otimes_k L$.
Thus $\varphi(h_f)$ and $\varphi(h_g)$ are elements of the free Lie algebra $Q\otimes_k L$ over $Q$ which are left dependent over $B$. By Lemma \ref{l2}, $\varphi(h_f)$ and $\varphi(h_g)$ are Lie dependent over $Q$. It is well-known (see, for example \cite{Shir1,Witt})  that in this case $\varphi(h_f)$ and $\varphi(h_g)$ are linearly dependent over $Q$, i.e., $\varphi(h_f)=\lambda \varphi(h_g)$. Consequently, $\mathrm{hdeg}(h_f)=\mathrm{hdeg}(h_g)>\mathrm{hdeg}(h_f-\lambda h_g)$.

Suppose that $h_f-\lambda h_g\neq 0$. Then $h_f-\lambda h_g$ and $h_g$ are left $B$-dependent again since so are $h_f$ and $h_g$. Notice that $h_f-\lambda h_g\in QL$. Consequently, repeating the same considerations as above we get
$\mathrm{hdeg}(h_f-\lambda h_g)=\mathrm{hdeg}(h_g)$, which is a contradiction.

Therefore $h_f=\lambda h_g$.
 Then,
\bes
\sum_{i\geq 1}(\partial_{e_i}(f)-\lambda  \partial_{e_i}(g)) h_{e_i}=0
\ees
and hence
\bes
\partial_{e_i}(f)-\lambda \partial_{e_i}(g)=0
\ees
for all $i\geq 1$. It is well-known (see, for example \cite{Umi24}) that in this case $f$ and $g$ are algebraically dependent.
$\Box$

Notice that the equivalence of conditions $(i)$ and $(iii)$ is proved in \cite{MLSh} and the equivalence of conditions $(iii)$ and $(v)$ for free Poisson algebras is proved in \cite{MakarU2}. All statements of this theorem for free Poisson algebras are established in \cite{UU10}.

\begin{theor}\label{t3} Every automorphism of the free Poisson field $P(x,y)$ in two variables $x,y$ over a field $k$ of characteristic $0$ is an automorphism of the field of rational functions $k(x,y)$, i.e.,
\bes
\mathrm{Aut}\,P(x,y) \cong \mathrm{Aut}\,k(x,y)=Cr_2(k).
\ees
\end{theor}
\Proof
 Let $\psi$ be an automorphism of $Q=P(x,y)$. Put $\psi(x)=f$ and $\psi(y)=g$.  We will show that  $f,g\in k(x,y)$.
 By Corollary \ref{c4}, it is sufficient to prove that $\mathrm{hdeg}(h_f)=\mathrm{hdeg}(h_g)=1$.

Suppose that $\mathrm{hdeg}(h_f)+\mathrm{hdeg}(h_g)\geq 3$. Note that $\Omega_Q=Q^e h_f+Q^e h_g=Q^e h_x+Q^e h_y$. Consequently, $Q^e h_f+Q^e h_g$ contains two elements $h_x$ and $h_y$ of $h$-degree $1$. It follows that
$\varphi(h_f)$ and $\varphi(h_g)$ are left dependent in the free associative algebra $Q\otimes_k \mathrm{Asso}\langle h_x,h_y\rangle$.

As in the proof of Theorem \ref{t2}, $\mathrm{hdeg}(h_f)=\mathrm{hdeg}(h_g)$ and there exist $0\neq \lambda\in Q$ such that $\mathrm{hdeg}(h_f-\lambda h_g)< \mathrm{hdeg}(h_f)$. Put $T=h_f-\lambda h_g$.
Then $\Omega_Q=Q^e T+Q^e h_g, \ \mathrm{hdeg}(h_f)+\mathrm{hdeg}(T)\geq 3$, and $T\in Q\otimes_k \mathrm{Lie}\langle h_x,h_y\rangle$. Hence $\mathrm{hdeg}(T)=\mathrm{hdeg}(h_g)$, which is impossible.
Therefore $\mathrm{hdeg}(h_f)=\mathrm{hdeg}(h_g) = 1$ and $f,g\in k(x,y)$. Similarly $\psi^{-1}(x),\psi^{-1}(y)\in k(x,y)$ and the restriction of $\psi$ on $k(x,y)$ is an automorphism.

The Theorem is proved since every automorphism of $k(x,y)$ can be uniquely extended to an automorphism of $P(x,y)$.
$\Box$

\bigskip

\begin{center}
{\bf\large Acknowledgments}
\end{center}

\hspace*{\parindent}

We are grateful to Max-Planck Institute f\"ur Mathematik for
hospitality and excellent working conditions, where part of this work has been done.

\end{document}